\documentclass[11pt]{article}
\usepackage[margin=1in]{geometry}
\usepackage{amsthm,amsmath}
\usepackage{latexsym,amssymb,amsmath}

\usepackage{tikz}
\usetikzlibrary{positioning,fit,backgrounds,calc}
\usepackage{xcolor}

\usepackage{graphicx,float,color,fancybox,shapepar,setspace,hyperref}
\usepackage{subfigure}
\usepackage{pgf,tikz}
\usepackage[affil-it]{authblk}
\usepackage{indentfirst}
\usepackage{hyperref}
\usepackage[section]{placeins}
\hypersetup
{
colorlinks=true,
linkcolor=blue,
filecolor=blue,
urlcolor=blue,
citecolor=cyan,
}
\usetikzlibrary{arrows}
\newtheorem{theorem}{Theorem}     
\newtheorem{lemma}{Lemma}         
\newtheorem{corollary}{Corollary} 

\newtheorem{conjecture}{Conjecture}

\newtheorem{problem}{Problem}

\newcounter{remarkcounter}

\def\~{\sim}

\begin{document}
\title{A note on tree-cycle Ramsey numbers}
\author{Ting HUANG$^1$, Yanbo ZHANG$^2$, Yaojun CHEN$^{1,}$\footnote{Corresponding author. Email: yaojunc@nju.edu.cn}\\
{\small $^1$School of Mathematics, Nanjing University, Nanjing 210093, China}\\
{\small $^2$School of Mathematical Sciences, Hebei Normal University, Shijiazhuang 050024, China}}
\date{}
\maketitle
	
\begin{quote}
	{\bf Abstract:}
    Let $R(T_n,C_m)$ denote the Ramsey number of a tree $T_n$ on $n$ vertices versus a cycle $C_m$ of length $m$. Burr, Erd\H{o}s, Faudree, Rousseau, and Schelp (1982) asked for the least function $f(m)$ such that $R(T_n,C_m)=2n-1$ for every odd $m\ge 3$ whenever $n\ge f(m)$. They proved that $f(m)\le 756m^{10}$. This bound was later improved to $25m$ by Brennan (2016) and to $4m-8$ by Fan and Lin (2025). In this note, we show that $f(m)\le 2m-4$ by using a different method and conjecture that $f(m)=\lceil (2m-1)/3\rceil$.


\end{quote}
\section{Introduction}
 Given two graphs $G$ and $H$, the \emph{Ramsey number} $R(G, H)$ is defined as the smallest positive integer $N$ such that, for every red-blue edge coloring of the complete graph $K_N$, there exists either a red subgraph isomorphic to $G$ or a blue subgraph isomorphic to $H$. For a graph $H$ with chromatic number $\chi(H)$, the \emph{chromatic surplus} $s(H)$ is defined to be the smallest size of a color class over all proper $\chi(H)$-colorings of $H$. For a connected graph $G$ with $n \geq s(H)$ vertices, Burr~\cite{Burr1981} established the following general lower bound for $R(G,H)$: 
\begin{equation} \label{lower}
        R(G, H)\ge (\chi(H)-1)(n-1)+s(H)\,.
\end{equation}

Burr~\cite{Burr1981} also defined $G$ to be  \emph{$H$-good} if equality holds in \eqref{lower}. Let $K_m$, $P_m$, and $C_m$ denote a complete graph, a path and a cycle on $m$ vertices, respectively, and let $T_m$ be an arbitrary tree on $m$ vertices. Prior to the formal definition of Ramsey goodness, Chv\'atal \cite{Chvatal1977} had established that all $n$-vertex trees $T_n$ satisfy the $K_m$-goodness property, as demonstrated in the following result.
 \begin{theorem}[Chv\'atal~\cite{Chvatal1977}]\label{Chv}  $R(T_n,K_m)=(n-1)(m-1)+1$.
 \end{theorem}

For the case $m=3$, Chv\'atal's result implies that $R(T_{n},K_{3}) = R(T_{n},C_{3})=2n-1$. Burr et al. extended this result by replacing $K_{3}$ with an odd cycle $C_{m}$ as follows.

\begin{theorem}[Burr, Erd\H{o}s, Faudree, Rousseau, and Schelp \cite{BEFRS1982}]
Let $T_n$ be an arbitrary tree on $n$ vertices. Then for odd $m\geq 3$ and $n\geq 756m^{10}$, $R(T_n,C_m)=2n-1$.
\end{theorem}

They also proposed the following problem.

\begin{problem}[Burr et al. \cite{BEFRS1982}] \label{problem}
What is the least function $f(m)$ such that $R(T_n,C_m)=2n-1$ for every odd $m\ge 3$ whenever $n\ge f(m)$?
\end{problem}

In 2016, Brennan showed that $n\geq 25m$ will suffice.

\begin{theorem}[Brennan \cite{Brennan2016}]
$R(T_{n},C_{m})=2n-1$ for odd $m\geq 3$ and $n\geq 25m$.
\end{theorem}

Recently, Fan and Lin improved Brennan's bound by reducing the requirement on $n$ from $n\geq 25m$ to $n\geq 4m-8$.

\begin{theorem}[Fan and Lin \cite{FanLin2025}]
$R(T_{n},C_{m})=2n-1$ for odd $m\geq 3$ and $n\geq 4m-8$.
\end{theorem}

In this note, by using a method different from the ones applied in \cite{Brennan2016,FanLin2025}, we further reduce the lower bound of $n$ to $2m-4$.
\begin{theorem} \label{mr}
$R(T_{n},C_{m})=2n-1$ for odd $m\geq 3$ and $n\geq 2m-4$.
\end{theorem}

For the case when $T_n$ is a path, Faudree, Lawrence, Parsons, and Schelp \cite{Faudree1974} showed that $R(P_{n},C_{m})=\max\{2n-1,m+\lfloor n/2 \rfloor-1\}$ for $m\geq n\geq 2$ and odd $m$. In particular, when $n\leq 2(m-1)/3$, the second term dominates, giving $R(P_{n},C_{m})=m+\lfloor n/2\rfloor-1>2n-1$, which implies $f(m)>2(m-1)/3$. Combining this with the upper bound $f(m)\leq 2m-4$ given in Theorem \ref{mr}, we have
\[
\lceil (2m-1)/3 \rceil\le f(m)\leq 2m-4.
\]

Furthermore, we conjecture that the lower bound is sharp for odd $m\ge 5$.

\begin{conjecture}
$R(T_{n},C_{m})=2n-1$ for odd $m\geq 5$ and $n\geq \lceil (2m-1)/3 \rceil$.

\end{conjecture}

\section{Preliminaries}

The symbols $|G|$, $e(G)$, $\delta(G)$, and $\Delta(G)$ denote the number of vertices, the number of edges, the minimum degree, and the maximum degree of the graph $G$, respectively. 
 By $\Delta_2(G)$, we
mean the second-largest degree. That is, if $(d_1,d_2,\ldots ,d_n)$ is the degree sequence  of $G$ with
$d_1\geq d_2 \geq \cdots \geq d_n$, then $\Delta_2(G)=d_2$. Let $K_{s,t}$ be the complete bipartite graph with parts of sizes $s$ and $t$.
For a red-blue edge-coloring of a $K_N$, we let $R$ and $B$ denote the spanning subgraphs of $K_N$, whose edge sets consist of all red and all blue edges of $K_N$, respectively. Given a subset $X\subseteq V(G)$, we write $G[X]$ for the subgraph of $G$ induced by $X$.
For a vertex $v\in V(G)$, $d_R(v)$ and $d_B(v)$ denote
the red degree and blue degree of $v$, that are defined as the number of red and blue edges incident to $v$, respectively. We denote by $N_{R}(v)$ (resp. $N_{B}(v)$) the set of vertices adjacent to $v$ by a red (resp. blue) edge.

Our proof also needs the following three important lemmas.
The first one is due to Dobson. From this lemma, we get that when a  graph contains no $K_{2,n}$, its complement graph could contain a large class of trees.
\begin{lemma}[Dobson \cite{Dobson2002}]\label{lem:noK2n}
    Let $n \geq 2$, and let $G$ be a graph of order $2n-1$ such that $\overline{G}$ contains no $K_{2,n}$. Then $G$
contains every tree $T$ of order $n$ such that $\Delta(T) \leq  \Delta(G)$ and $\Delta_2(T) \leq \delta(G)$.
\end{lemma}

The second one is due to H\"aggkvist, which gives an upper bound for $R(T_n, P_m)$.

\begin{lemma}[H\"aggkvist \cite{Haggkvist1989}]
    \[
R(K_{n_1,n_2}, P_m) \leq n_1 + n_2 + m- 2.
\]
\end{lemma}
    Since every tree is a subgraph of a complete bipartite graph, the following holds.
    \begin{corollary} \label{TP}
    For $m\geq 2$ and $n\geq 2$,
\[
R(T_n, P_m) \leq n + m- 2.
\]
\end{corollary}

 The concept of pancyclicity for graphs was introduced and studied by Bondy \cite{Bondy1971} in the early 1970s. A graph $ G $ is pancyclic if it contains cycles of every length between 3 and $|G|$. Bondy established several sufficient conditions for a graph to be pancyclic. A typical degree condition is given below.

\begin{lemma}[Bondy \cite{Bondy1971}] \label{pancyclic}
Let $G$ be a graph of order $n$ with $\delta(G) \geq n/2 $. Then $ G $ is pancyclic, or $ G = K_{r,r} $ with $ r = n/2 $.
\end{lemma}

\section{Proof of Theorem \ref{mr}}

From \eqref{lower}, we have $R(T_n, C_m) \geq 2n - 1$ for every odd integer $m \geq 3$. In what follows, we show that if $m \geq 3$ is odd and $n \geq 2m-4$, then $R(T_n, C_m) \leq 2n - 1$. Let $N=2n-1$. Assume that $K_N$ has a red-blue edge coloring such that it contains neither a red $T_n$ nor a blue $C_m$. 

If $\delta(R) \geq n - 1$, then $T_n$ can be embedded into $R$ clearly. Hence $\delta(R) \leq n - 2$, which implies $\Delta(B) = N- 1 - \delta(R) \geq n$. Let $u \in V(K_N)$ be such that $d_B(u) = \Delta(B)$. Then there must be one edge in $B[N_B(u)]$, which together with $u$ forms a blue $C_3$, a contradiction. Hence we may assume 
$m\geq 5$.

Then we claim that $\Delta(R)\ge n-1$ and $\delta(R)\ge n/2$.

If $\delta(B)\ge n>N/2$, then by Lemma~\ref{pancyclic}, since $N=2n-1$ is odd,
$B$ is pancyclic, and hence contains a blue $C_m$, a contradiction. Thus
$\delta(B)\le n-1$, and so \begin{equation}
\Delta(R) \geq N - 1 - (n - 1) = n - 1 . \label{eq:delta_bound}
\end{equation}

Suppose next that $\delta(R)<n/2$. Let $u$ be a vertex with
$d_R(u)=\delta(R)$. Since $d_R(u)$ is an integer,
$d_R(u)\le \left\lceil n/2\right\rceil-1$.
Therefore
$d_B(u)\ge 2n-1-\left\lceil n/2\right\rceil
= n-1+\left\lfloor n/2\right\rfloor$.
As $n\ge 2m-4$, we have $\lfloor n/2\rfloor\ge m-2$, and hence
$d_B(u)\ge n+m-3$.
By Corollary~\ref{TP},
$R(T_n,P_{m-1})\le n+(m-1)-2=n+m-3$.
Thus $B[N_B(u)]$ contains a blue $P_{m-1}$, unless $R[N_B(u)]$ contains
a red $T_n$. The latter is impossible, while the former together with $u$
forms a blue $C_m$, a contradiction. Hence \begin{equation}
\delta(R) \geq \frac{n}{2} . \label{eq:delta2_bound}
\end{equation}

Notably, by $e(T_n)=n-1$, we have $\Delta_2(T_n)\leq n/2$ and $\Delta(T_n)\leq n-1$ for all $n$-vertex trees $T_n$. By 
(\ref{eq:delta_bound}) and (\ref{eq:delta2_bound}), we have $\Delta_2(T_n)\leq \delta(R)$ and $\Delta(T_n)\leq \Delta(R)$ for all $n$-vertex trees $T_n$. If $K_N$ contains no blue $K_{2,n}$, then
applying Lemma \ref{lem:noK2n} on $G=R$ and $\overline{G}=B$, we conclude that $R$ contains every tree $T_n$, that is, $K_N$ contains a red copy of $T_n$, a contradiction.
Hence, $K_N$ contains a blue $K_{2,n}$.

Let the two partite sets of the blue $K_{2,n}$ be $\{u_1,u_2\}$ and $U$, where $|U|=n$. Set $V=V(K_N)$ and
$W = V \setminus (U \cup \{u_1, u_2\})$. Then $|W|=n-3$. 

Clearly $B[U]$ contains no $P_{m-3}$. Indeed, if
$z_1z_2\cdots z_{m-3}$ is a blue path in $U$, then since
$n-(m-3)=n-m+3>0$, we may choose $y\in U\setminus\{z_1,\ldots,z_{m-3}\}$.
Then
$u_1z_1z_2\cdots z_{m-3}u_2yu_1$
is a blue $C_m$, a contradiction.

Let $T^1$ be a subtree obtained from $T_n$ by repeatedly deleting leaves
until exactly $n-m+5$ vertices remain. By Corollary~\ref{TP},
\[
R(T_{n-m+5},P_{m-3})\le (n-m+5)+(m-3)-2=n=|U|.
\]
Since $B[U]$ contains no $P_{m-3}$, the graph $R[U]$ contains a red copy
of $T^1$. If $m=5$, then $T^1=T_n$, a contradiction. Thus we may assume
$m\ge 7$.

Starting from this red copy of $T^1$, greedily extend the embedding along
a leaf-adding order from $T^1$ to $T_n$ in $R$. Since $R$ contains no red $T_n$,
the procedure must fail. Hence, for the current embedded subtree $H$,
there exists a vertex $u_3\in V(H)$ such that $u_3$ is blue-adjacent to
every vertex in $V(K_N)\setminus V(H)$. Moreover, $|V(H)|\le n-1$.

Since the initial copy of $T^1$ lies entirely in $U$, and at most
$(n-1)-(n-m+5)=m-6$
vertices have been added before the failure, we have
$|W\cap V(H)|\le m-6$.

Set $A=N_B(u_3)\cap W$.
Since $W\setminus V(H)\subseteq A$, we get
\[
|A|\ge |W|-|W\cap V(H)|
\ge (n-3)-(m-6)=n-m+3.
\]

We claim that $u_3\notin\{u_1,u_2\}$. Otherwise, since $u_3$ is blue-adjacent
to every vertex of $U$ and $A$, we have
$d_B(u_3)\ge |U|+|A|\ge n+(n-m+3)=2n-m+3$.
Consequently
\[
d_R(u_3)\le 2n-2-(2n-m+3)=m-5<n/2,
\]
where the last inequality follows from $n\ge 2m-4$. This contradicts
$\delta(R)\ge n/2$. Hence $u_3\notin\{u_1,u_2\}$.

Moreover, since $|V(H)|\le n-1$ and $|U|=n$, there exists
$v_1\in U\setminus V(H)$.
By the choice of $u_3$, the edge $u_3v_1$ is blue. Define
$$U_0=U\setminus\bigl(\{v_1\}\cup(\{u_3\}\cap U)\bigr).$$
Then $|U_0|\ge n-2$. Since $U_0\cap A=\emptyset$ and $n\ge 2m-4$, we have
$$|U_0\cup A|\ge (n-2)+(n-m+3)=2n-m+1\ge n+m-3.$$
Again by Corollary~\ref{TP},
$R(T_n,P_{m-1})\le n+m-3$, this gives
$$|U_0\cup A| \ge R(T_n,P_{m-1}).$$
Since $R$ contains no red $T_n$, the graph $B[U_0\cup A]$ contains a blue
path $P_{m-1}$, say
$P=z_1z_2\cdots z_{m-1}$.

We now derive a blue $C_m$ according to the positions of the endpoints of
$P$.

If $z_1,z_{m-1}\in U_0$, then
$u_1z_1z_2\cdots z_{m-1}u_1$
is a blue $C_m$.

If $z_1,z_{m-1}\in A$, then
$u_3z_1z_2\cdots z_{m-1}u_3$
is a blue $C_m$.

If one of $z_1,z_{m-1}$ lies in $U_0$ and the other in $A$,
 without loss of generality, let
$z_1\in U_0,~ z_{m-1}\in A$.
Consider the position of  $z_{m-3}$.

If $z_{m-3}\in U_0$, then
$z_1z_2\cdots z_{m-3}u_1v_1u_2z_1$
is a blue $C_m$.
 
 If $z_{m-3}\in A$, then
$z_1z_2\cdots z_{m-3}u_3v_1u_1z_1$
is a blue $C_m$.

In all cases we obtain a blue $C_m$, a contradiction. 

The proof of Theorem \ref{mr} is complete.

\section*{Acknowledgement}
\noindent This research is supported by National Key R\&D Program of China under grant number 2024YFA1013900 and NSFC under grant number 12471327.
\section*{Declaration}
	
\noindent$\textbf{Conflict~of~interest}$
The authors declare that they have no known competing financial interests or personal relationships that could have appeared to influence the work reported in this paper.
\vskip 2mm	
\noindent$\textbf{Data~availability}$
No data was used for the research described in the article.

\end{document}